\documentclass[12pt]{amsart}
\usepackage[margin=1.0in]{geometry}
\usepackage{ifxetex,ifluatex}
\newif\ifxetexorluatex
\ifxetex
  \xetexorluatextrue
\else
  \ifluatex
    \xetexorluatextrue
  \else
    \xetexorluatexfalse
  \fi
\fi

\ifxetexorluatex
  \usepackage{fontspec}
\else
  \usepackage[T1]{fontenc}
  \usepackage[utf8]{inputenc}
  \usepackage{lmodern}
\fi
\usepackage{hyperref, amsmath, amssymb, mathtools}
\usepackage{amsthm}
\usepackage{amsfonts}
\usepackage{amsmath}
\usepackage{comment}
\usepackage{chngcntr}
\newtheorem{thm}{Theorem}
\newtheorem{lemma}[thm]{Lemma}
\newtheorem{corollary}[thm]{Corollary}

\theoremstyle{definition}
\newtheorem{example}{Example}

\newtheorem*{remark}{Remark}
\numberwithin{thm}{section}
\DeclareMathOperator{\Tr}{Tr}
\DeclareMathOperator{\Frob}{Frob}
\title{Generating Functions for Power Moments of Elliptic Curves over $\mathbb{F}_{p}$}

\author[K. Gallagher]{Katherine Gallagher}
\address{Department of Mathematics, University of Notre Dame, Notre Dame, IN  46556}
\email{kgalla17@nd.edu}

\author[L. Li]{Lucia Li}
\address{Department of Mathematics, Wellesley College, Wellesley, MA 02481}
\email{lucia.li@wellesley.edu}

\author[N. Sweeting]{Naomi Sweeting}
\address{Department of Mathematics, University of Chicago, Chicago, IL 60637}
\email{nsweeting@uchicago.edu}

\author[K. Vassilev]{Katja Vassilev}
\address{Department of Mathematics, Princeton University, Princeton, NJ 08544}
\email{kdv@princeton.edu}

\author[K. Woo]{Katharine Woo}
\address{Department of Mathematics, Stanford University, Stanford, CA 94305}
\email{katywoo@stanford.edu}

\begin{abstract}
Seminal works by Birch and Ihara gave formulas for the $m$th power moments of the traces of Frobenius endomorphisms of elliptic curves over $\mathbb{F}_{p}$ for primes $p \geq 5$. Recent works by Kaplan and Petrow generalized these results to the setting of elliptic curves that contain a subgroup isomorphic to a fixed finite abelian group $A$. We revisit these formulas and determine a simple expression for the zeta function $Z_p(A; t)$, the generating function for these $m$th power moments. In particular, we find that
\[
Z_p(A;t) = \frac{\widehat{Z}_p(A; t)}{\hspace*{-0.3cm}\displaystyle \prod_{a \in \Frob_p(A)}\hspace*{-0.3cm}(1 - at)},\]
where $\Frob_p(A) \coloneqq \{ a \, \colon -2\sqrt{p} \leq a \leq 2\sqrt{p}\, \text{ and }  a \equiv p+1 \pmod{|A|}\}$, and $\widehat{Z}_p(A;t)$ is an easily computed polynomial that is determined by the first $\Big\lceil\frac{2\lfloor 2\sqrt{p}\rfloor}{|A|}\Big\rceil$ power moments.
These rational zeta functions have two natural applications. We find rational generating functions in weight aspect for traces of Hecke operators on $S_k(\Gamma)$ for various congruence subgroups $\Gamma$. We also prove congruence relations for power moments by making use of known congruences for traces of Hecke operators.

\end{abstract}
 
\section{Introduction and Statement of Results} \label{sec1}
For an elliptic curve in Weierstrass normal form $E(\alpha, \beta) \colon y^2z = x^3 - \alpha xz^2 - \beta z^3$ and a prime $p \geq 5$, we count the $\mathbb{F}_p$-rational points $N_E(p)$ by \begin{equation}\numberwithin{equation}{section}
a_E(p) \coloneqq p + 1 - N_E(p) = -\sum_{x = 0}^{p-1} \left(\frac{x^3 - \alpha x - \beta}{p}\right),
\end{equation} where $\big(\frac{\cdot}{p}\big)$ denotes the Legendre symbol. For an elliptic curve $E$ and a prime $p$, $a_{E}(p)$ is the trace of the Frobenius endomorphism, and its magnitude is bounded by $2\sqrt{p}$ according to Hasse's theorem. We can write $a_{E}(p)=2 \sqrt{p} \cos(\theta_p)$ for a unique $0\leq \theta_p\leq \pi$. For a fixed $E$, the distribution of the traces of Frobenius endomorphisms over varying prime fields is the content of the famous Sato-Tate conjecture, which was  proved by Taylor and his collaborators in \cite{STI, STII, STIII}.  Namely, for $E$ without complex multiplication and any $0 \leq c < d \leq \pi$, their work implies that the density of primes such that $c \leq \theta_p \leq d$ is given by \[\frac{2}{\pi} \int_c^d \sin^2 t \ dt.\]

Prior to this work, Birch \cite{Birch} and Ihara \cite{Ihara} investigated the distribution of the $a_E(p)$ by instead fixing a finite field $\mathbb{F}_p$ and varying the elliptic curve. To obtain this distribution, they considered the $m$th power moments\footnote{The reader should be aware that our sum, unlike its analogue in Birch's paper, excludes singular curves. This change does not affect the asymptotics.} defined for $p \geq 5$ and $m \geq 0$ by 

\begin{equation*}
M_p(m) \coloneqq \sum_{\substack{\alpha,\beta \in \mathbb{F}_{p} \\ 4\alpha^3 - 27\beta^2 \neq 0}} \bigg[\sum_{x \in \mathbb{F}_p}\left(\frac{x^3 - \alpha x - \beta}{p}\right)\bigg]^m.
\end{equation*}

Birch's main result, reformulated explicitly in terms of the Catalan numbers $C_m$, gives an asymptotic for even moments as $m\rightarrow \infty$ (all odd moments are zero):
\begin{equation} \label{asymptotics}
M_p(2m) \sim C_m \cdot p ^{m + 2} = \frac{(2m) !}{m! (m + 1)!}\cdot p^{m + 2}.
\end{equation}
This result has recently been generalized by Kaplan and Petrow \cite{KaplanPetrow}, and in a slightly different guise by Kowalski \cite{Kowalski}, to compute power moments for only those elliptic curves whose group of $\mathbb{F}_p$-rational points contains a subgroup isomorphic to a fixed abelian group $A$. These moments are defined by \begin{equation} M_p(A; m) \coloneqq \hspace*{-.3cm} \sum_{\substack{\alpha, \beta \in \mathbb{F}_p \\ 4 \alpha^3 - 27 \beta^2 \neq 0 \\ A \hookrightarrow E(\alpha, \beta)}} \bigg[\sum_{x \in \mathbb{F}_p}\left(\frac{x^3 - \alpha x - \beta}{p}\right)\bigg]^m.
\end{equation} For example, $M_p(A; 0)$ is the number of $(\alpha, \beta) \in \mathbb{F}_p^2$ so that $y^2z = x^3 - \alpha xz^2 - \beta z^3$ is a nonsingular elliptic curve whose group of $\mathbb{F}_p$-rational points contains a subgroup isomorphic to $A$.

We can also express the moments as \begin{equation} \label{momentdef}M_p(A; m) = (p^2-p) \mathbb{E}_p \big(\Phi_A a_{E}(p)^m\big),\end{equation} where $\mathbb{E}_p$ denotes the expectation over all $(\alpha, \beta) \in \mathbb{F}_p^2$ such that $4\alpha^3 - 27\beta^2 \neq 0$ and $\Phi_A$ denotes the indicator function of $A \hookrightarrow E(\alpha,\beta)$. We study the generating function of these moments, the zeta function \begin{equation}
Z_p(A; t) \coloneqq \sum_{m =0}^\infty M_p(A; m) t ^ m. 
\end{equation}

We prove that $Z_p(A; t)$ is a simple rational function in $t$ that can be calculated from the first few values of $M_p(A; m)$. To be precise, for a prime $p \geq 5$, we define \begin{equation}\label{delta}\delta_A(p)\coloneqq \#\left\{a \, \colon  0 < |a| \leq 2\sqrt{p}\, \textrm{ and } a \equiv p+1 \pmod{|A|}\right\},\end{equation} and we observe that: \[\delta_A(p) \leq \bigg \lceil\frac{2\lfloor 2\sqrt{p}\rfloor}{|A|} \bigg \rceil.\] We also define \[ \Frob_p(A) \coloneqq \left\{a \,\colon  -2\sqrt{p} \leq a \leq 2\sqrt{p}\, \textrm{ and } \, a \equiv p+1 \pmod{|A|} \right\}. \] Here, $\Frob_p(A)$ should be thought of as the set of potential traces of Frobenius endomorphisms of elliptic curves $E$ whose group of $\mathbb{F}_p$-rational points contains a subgroup isomorphic to $A$. Using the $m$th power moments $M_p(A;0), \ldots , M_p(A; \delta_A(p))$, we define a sequence of integers $c_p(A;n)$ by \begin{equation*}\sum_{n = 0}^\infty c_p(A;n)t^n \coloneqq \left(M_p(A;0) + M_p(A;1)t + \cdots + M_p(A; \delta_A(p)) t^{\delta_A(p)} \right) \hspace*{-0.2cm}\displaystyle\prod_{a \in \Frob_p(A)} \hspace*{-.3cm}(1 - at).\end{equation*} This sequence provides the coefficients of the polynomial \begin{equation}\widehat{Z}_p(A; t) \coloneqq c_p(A;0) + c_p(A;1)t + \cdots + c_p(A;\delta_A(p))t^{\delta_A(p)}. \end{equation} 

We now state our first theorem, which gives a simple form for the zeta function $Z_p(A; t)$. 
\begin{thm}\label{TheoremZetaGF}
If $p \geq 5$ is prime, then we have that \[Z_p(A;t) = \frac{\widehat{Z}_p(A; t)}{\displaystyle\prod_{a \in \Frob_p(A)}\hspace*{-.3cm}(1 - at)}.\]
\end{thm}

\begin{remark} Since $Z_p(A; t)$ is a rational function with denominator of degree $\delta_A(p)$, the moments $M_p(A; m)$ satisfy a recurrence relation of length $\delta_A(p)$. In the special cases that $|A| \in \left\{1, 2\right\}$, all the odd moments are trivial, so the even moments satisfy a recurrence of length $\frac{\delta_A(p)}{2}$. In general, $\delta_A(p)$ can be thought of as a measure of the complexity of the power moment problem. 
\end{remark} 
\begin{remark}
It is well-known that only some $A$ can be isomorphic to subgroups of elliptic curves over $\mathbb{F}_p$. Namely, $A$ must be of the form $\mathbb{Z}/n_1\mathbb{Z} \times \mathbb{Z}/n_2\mathbb{Z}$ where $n_2 \mid n_1$.  Although Theorem \ref{TheoremZetaGF} is unnecessary if $A$ is not of this form, it still applies, as all the low-order moments $M_p(A; m)$ vanish, and so we have that $\widehat{Z}_p(A; t) = 0$. 
\end{remark}

In contrast to our generating function approach, the works by Birch, Ihara, Kaplan, and Petrow \cite{Birch, Ihara, KaplanPetrow} on power moments used versions of the Eichler-Selberg trace formula for the action of Hecke operators on certain spaces of modular forms. The connection between Hecke operators and power moments appears because, as Deuring \cite{Deuring} proved, any elliptic curve over $\mathbb{F}_p$ is the reduction of a curve with complex multiplication. The theory of complex multiplication in turn provides the link between counting elliptic curves and the Hurwitz class numbers which appear in the trace formula \cite{Cox}.

Thanks to this connection, Theorem \ref{TheoremZetaGF} can be used to compute the traces of Hecke operators. Our zeta functions $Z_p(A; t)$ allow us to explicitly determine the rational generating functions in weight aspect for all traces of Hecke operators $T_p$ on $S_k(\Gamma)$ for certain congruence subgroups $\Gamma$. We give a special case to simplify the calculations, but the machinery of \cite{KaplanPetrow} works in greater generality, as do our generating function methods. Let $A = \mathbb{Z}/n_1\mathbb{Z} \times \mathbb{Z}/n_2\mathbb{Z}$ with $n_2 \mid n_1$. The constants $c(A)$ and $d_p(A; k)$, and the congruence subgroup $\Gamma(n_1, n_2)$, are defined in \S \ref{Notation}. In particular, it turns out that $\Gamma(N, N) \simeq \Gamma(N)$ and $\Gamma(N, 1) = \Gamma_1(N)$.

\begin{thm}\label{TheoremTraceGF} 
If $p \geq 5$ is prime, $\gcd(p, |A|) = 1$, $p\equiv 1 \pmod{n_2}$, and $\gcd(p-1, n_1) = n_2$, then the generating function in weight aspect of traces of $T_p$ acting on $S_k(\Gamma(n_1,n_2))$ is given by: \begin{align*}\sum_{k = 2}^\infty \Tr\big(T_p | S_k \big(\Gamma(n_1, n_2)\big)\big)t^k = (p+1)t^2 &+ \frac{d_p(A; 0) t^2}{1 - t^2} + \frac{d_p(A; 1) t^3}{1 - t^2} \\ &- \frac{1}{c(A) }\cdot \frac{t^2}{(p-1)(1 + pt^2)} \cdot Z_p\left(A; \frac{t}{1 + pt^2}\right).\end{align*}
\end{thm}

\begin{corollary}
In the special case $A = \left\{0\right\}$, we obtain: \[\sum_{k=1}^{\infty} \Tr\big(T_{p}|S_{2k} \big(SL_2(\mathbb{Z})\big)\big)t^{2k} = (p+1)t^2 + \frac{-t^2}{1-t^{2}} - \frac{t^{2}}{(p-1)(1+pt^{2})} \cdot Z_{p} \left(\{0\}; \frac{t}{1+pt^2}\right).\]
\end{corollary}

\medskip
Our simple form for the generating function $Z_p(A;t)$ motivates us to investigate congruences in ${M_p(2m) := M_p(\left\{0\right\}; 2m)},$ whose generating function is denoted ${Z_p(t) := Z_p(\left\{0\right\}; t)}$. In particular, we prove the following theorem.
\begin{thm}\label{congruencethm} The following congruence relations are true:

\noindent (1) If $\ell\in \{2,3,5\}$, $p \geq 5$ is a prime such that $p \equiv -1 \pmod{\ell}$, and $m>0$, then \[\frac{M_p(2m)}{p-1} \equiv 0 \pmod{\ell}.\]
(2) If $p \geq 5$ is prime, $\mu \in \left\{4, 6, 8, 10, 14\right\}$, $m \equiv \mu \pmod{p-1}$, and $m>0$, then \[M_p(m) \equiv 1 \pmod{p}.\]
\end{thm}

The plan of the paper is as follows. First, \S \ref{sec2} introduces necessary background and notation. Next, \S \ref{sec3} provides proofs of the main theorems. Finally, \S \ref{sec4} gives examples of our results.
 
\section{Preliminaries} \label{sec2}

\subsection{Notation}\label{Notation}
First, we define some notation that appears in the statement of Theorem \ref{TheoremTraceGF}, as well as in \cite{KaplanPetrow}.  Let \[\psi(n) \coloneqq n\prod_{p\mid n} \left(1 + \frac{1}{p}\right),\] and let $\varphi(n)$ be the Euler totient function. Write $\delta_n(a, b)$ for the indicator function of ${a \equiv b \pmod{n}}$ and $\delta(a,b)$ for the indicator function of $a = b$. For $A = \mathbb{Z}/n_1 \mathbb{Z} \times \mathbb{Z}/n_2 \mathbb{Z}$ with $n_2 \mid n_1$, define \[c(A)\coloneqq \frac{\psi(n_1^2/n_2^2)}{\psi(n_1^2)\varphi(n_1)},\] and define $d_p(A; k)$ to be
\[\frac{-\varphi(n_1)}{4} \cdot \Big(1 + (-1)^k \delta_{n_1}(-1, 1) + \delta_{n_1} (p, 1) + (-1)^k \delta_{n_1} (p, -1)\Big) \cdot \hspace{-0.6cm}\sum_{\substack{\tau \mid n_1 n_2 \\ \gcd(\tau, n_1n_2/\tau)\mid n_2}}\hspace*{-0.7cm} \varphi\left(\gcd\left(\tau, \frac{n_1n_2}{\tau}\right)\right). \]

\subsection{Hecke operators and the Eichler-Selberg trace formula}
Note that the interested reader can refer to \cite{Ono} for more information on modular forms and Hecke operators. Let $M_{k}(\Gamma)$ (resp. $S_{k}(\Gamma)$) denote the complex vector space of modular forms (resp. cusp forms) of weight $k$ on a congruence subgroup $\Gamma \leq SL_{2}(\mathbb{Z})$. The Hecke operators $T_{p}$ are linear transformations that act on these spaces of modular forms. For example, if $\displaystyle f(z)= \sum_{n=0}^{\infty}a(n)q^{n} \in M_{k} (SL_2(\mathbb{Z}))$, we have: $$\displaystyle f | T_{p} \coloneqq \sum_{n=0}^{\infty}\big(a(pn)+p^{k-1}a(n/p)\big)q^{n},$$
where we set $a(n/p)=0$ if $p \nmid n$. For the action of $T_p$ on $M_k(\Gamma)$ for various congruence subgroups $\Gamma$, see \cite{DS}. Let $\Tr \big(T_{p}| S_k(\Gamma)\big)$ denote the trace of the Hecke operator $T_{p}$ acting on $S_{k}(\Gamma)$. We sometimes write $\Tr_k(T_p) \coloneqq \Tr \big(T_{p}| S_k\big(SL_2(\mathbb{Z})\big)\big)$. Of particular interest to us is the congruence subgroup $\Gamma(N, M)$ for $M \mid N$ defined as follows:
\begin{equation}
\Gamma(N, M)\coloneqq \left\{\begin{pmatrix} a & b \\ c & d \end{pmatrix} \in SL_2(\mathbb{Z})  \colon a \equiv d \equiv 1 \pmod{N}, \, c \equiv 0 \pmod{NM} \right\}.
\end{equation}

The Eichler-Selberg trace formula is central to previous works on moments. We recall the normalized Chebyshev polynomials, \begin{equation}\label{Chebyshev}
U_{k-2}(a,b) \coloneqq \frac{\xi^{k-1} - \overline{\xi}^{k-1}}{\xi-\overline{\xi}},
\end{equation}
where $\xi$ is a solution to the polynomial $X^2-a X+ b=0$. The Eichler-Selberg trace formula for a prime Hecke operator $T_p$ on cusp forms on the full modular group is

\begin{equation*}
1 + \Tr_{2k}(T_p) = -\frac{1}{2} \sum_{|a| \leq 2 \sqrt{p}}  H(a^2 - 4 p) U_{2k -2},
\end{equation*} where $H(\cdot)$ denotes the Hurwitz class number. From this, Birch \cite{Birch} derived a simple relationship between traces and moments, valid for $k > 1$:
\begin{equation*} \label{Bircheq}
\Tr_{2k}(T_p) = -1 - \sum_{m = 0}^{k - 1} \binom{k + m - 1}{2m} \frac{M_p(2m)}{p - 1} (- p)^{k - 1 - m}.
\end{equation*}

For traces of Hecke operators over $\Gamma (n_1, n_2)$,  Kaplan and Petrow \cite{KaplanPetrow} proved a  more refined version of the Eichler-Selberg trace formula. Using this formula, they compute the expectations of traces of Frobenius endomorphisms weighted by Chebyshev polynomials, $\mathbb{E}_p\big(\Phi_A U_{k-2}(a_E,p)\big)$, in terms of the traces $\Tr\big(T_p | S_k\big(\Gamma(n_1,n_2)\big)\big)$. When $p \geq 5$ is prime and under the other assumptions of Theorem \ref{TheoremTraceGF}, Kaplan and Petrow's main result reduces to: 

\begin{equation} \label{KaplanPetroweq}
    \Tr\big(T_p | S_k\big(\Gamma(n_1,n_2)\big)\big) = (p+1)  \delta(k,2) - \frac{p\mathbb{E}_p\big(\Phi_A U_{k-2}(a_E,p)\big)}{c(A)} + d_p(A;k).
\end{equation}

\section{Proofs} \label{sec3}
\subsection*{Proof of Theorem \ref{TheoremZetaGF}}
By definition, we have that \\
\begin{equation*}
M_p(A;m) = \sum_{|a| \leq 2\sqrt{p}} a^{m} S_a(A),
\end{equation*}
where $S_{a}(A)\coloneqq \# \left\{ (\alpha, \beta) \in \mathbb{F}_p^2 \, \colon  A \hookrightarrow E(\alpha, \beta) 
\text{ and } a_{E}(p) = a\right\}$. Using the definition of $Z_p(A;t)$, we have:
\begin{equation*}
Z_p(A;t) = \sum_{m=0}^\infty M_p(A;m) t^m = \sum_{m=0}^{\infty} \sum_{|a|\leq 2\sqrt{p}}(at)^m S_a(A) = \sum_{|a|\leq 2\sqrt{p}} \frac{S_a(A)}{1-at}.
\end{equation*}
Because $p+1-a_E(p)= N_E(p)$, the condition that $A \hookrightarrow E(\alpha, \beta)$ implies that $a \equiv p+1 \pmod{|A|}$. Of course, this condition is necessary but not sufficient. Hence, $S_a(A) = 0$ whenever $a \not\equiv p+1 \pmod{|A|}$. So, our expression becomes:
\begin{equation*}
 Z_p(A;t) = \sum_{a \in \Frob_{p}(A)}\frac{S_a(A)}{1-at}.
\end{equation*}
Note that the denominator is a polynomial of degree $\delta_A(p)$, as defined in (\ref{delta}). 
In addition, the degree of the numerator is less than or equal to $\delta_A(p)$. Thus, if we define $$\sum_{n = 0}^\infty c_{p}(A;n)t^n \coloneqq \left(M_p(A;0) + M_p(A;1)t + \cdots + M_p(A; \delta_A(p)) t^{\delta_A(p)} \right) \displaystyle{\prod_{ a \in \Frob_{p}(A)}(1 - at)},$$ 
then $\widehat{Z}_p(A; t) \coloneqq c_p(A;0) + c_p(A;1)t + \dots + c_p(A;\delta_A(p))t^{\delta_A(p)}$ is the numerator of the rational function $Z_p(A;t)$.
  \qed
\subsection*{Proof of Theorem \ref{TheoremTraceGF}} For this proof, we  fix a prime $p \geq 5$ and denote $a_E \coloneqq a_E(p)$ for simplicity. Let $\mathbb{E}_p\big(\Phi_A U_{k-2}(a_E,p)\big)$ be the expectation of the normalized Chebyshev polynomials $U_{k-2}(a_E,p)$ defined in (\ref{Chebyshev}).
Making use of (\ref{KaplanPetroweq}) from \cite{KaplanPetrow}, we compute the generating function:
\begin{align*}
    \sum_{k=2}^\infty \Tr \big(T_p | S_k\big(\Gamma(n_1,n_2)\big)\big) t^k &= \sum_{k=2}^\infty \Big((p+1)\delta(k,2) -\frac{p\mathbb{E}_p\big(\Phi_A U_{k-2}(a_E,p)\big)}{c(A)} + d_p(A;k)\Big)t^k \\
    &= (p+1)t^2 + \sum_{k=2}^\infty d_p(A;k)t^k - \frac{p}{c(A)} \sum_{k=2}^\infty \mathbb{E}_p\big(\Phi_A U_{k-2}(a_E,p)\big)t^k.
\end{align*}
By linearity of expectation, we rewrite part of this expression as: \begin{align*}\sum_{k = 2}^\infty \mathbb{E}_p\big(\Phi_A U_{k-2}(a_E,p)\big)t^k &= \mathbb{E}_p \Big(\Phi_A \sum_{k=2}^\infty U_{k-2}(a_E,p)t^k\Big) \\ & = \mathbb{E}_p \Big( \frac{\Phi_A t^2}{pt^2 - a_E t + 1}\Big)\\
    &= \frac{t^2}{(1+pt^2)} \cdot \mathbb{E}_p\left(\Phi_A\sum_{\ell=0}^\infty \left(\frac{a_Et}{1+pt^2}\right)^\ell \right). \\
    \end{align*}
At this point, we rewrite this in terms of moments using (\ref{momentdef}). Indeed, we have that
\[   \mathbb{E}_p\left(\Phi_A\sum_{\ell=0}^\infty \left(\frac{a_Et}{1+pt^2}\right)^\ell \right)  = \sum_{\ell=0}^\infty \frac{M_p(A;\ell)}{p(p-1)} \left(\frac{t}{1+pt^2}\right)^\ell = \frac{1}{p(p-1)}\cdot  Z_p\left(A; \frac{t}{1+pt^2}\right).\] Finally, we notice that $d_p(A;k)$ depends only on the parity of $k$. Putting everything together, we obtain:
\begin{align*}
    \sum_{k=2}^\infty \Tr \big(T_p|S_k\big(\Gamma(n_1,n_2)\big)\big)t^k = (p+1)t^2 &+ \frac{d_p(A;0)t^2}{1-t^2} + \frac{d_p(A;1)t^3}{1-t^2}  \\
    &-\frac{t^2}{(p-1)c(A)(1+pt^2)} \cdot Z_p \left( A; \frac{t}{1+pt^2} \right).
\end{align*}
\qed
\medskip

\subsection*{Proof of Theorem \ref{congruencethm}}
To prove Theorem \ref{congruencethm}, we define a natural decomposition of $Z_{p}(t)$. We then prove a series of lemmas about the decomposition in order to deduce the statement of the theorem.

Birch \cite{Birch} gave the following formula:
\begin{equation*}\label{BirchMp}
\frac{M_p(2m)}{p-1} = \frac{(2m)!}{m!(m+1)!}p^{m+1} - \sum_{k = 1}^m \frac{(2k+1)(2m)!}{(m-k)!(m+k+1)!} p^{m-k}\big(1 + \Tr_{2k + 2}(T_p)\big).
\end{equation*}
Motivated by this formula, we define what we term the combinatorial  and the Hecke part of the moments: 

\begin{equation*}
M_p^{{C}}(2m) \coloneqq (p-1)\left(\frac{(2m)!}{m!(m+1)!}p^{m+1} - \sum_{k = 1}^m \frac{(2k+1)(2m)!}{(m-k)!(m+k+1)!} p^{m-k}\right), \label{Mppoly}
\end{equation*}
\begin{equation*}
M_p^{H}(2m)\coloneqq -(p-1)\sum_{k = 1}^m \frac{(2k+1)(2m)!}{(m-k)!(m+k+1)!}p^{m-k} \Tr_{2k+2} (T_p). \label{Mphecke}
\end{equation*} 
For odd moments $2m + 1$, we set  $M_p^{C}(2m+1) = M_p^{H} (2m+1) = 0$. The corresponding generating functions are defined as
\begin{align*}
Z_p^{C}(t) \coloneqq \sum_{m = 0}^\infty M_p^{C} (m)t^m 
\,\,\,\,\text{  and  }\,\,\,\,\,
Z_p^{H}(t) \coloneqq \sum_{m=0}^\infty M_p^{H}(m)t^m.
\end{align*}
Clearly, $Z_p(t)$ decomposes as $Z_p^C(t) + Z_p^H(t)$. 

Our first lemma gives a closed formula for  $Z_p^C(t)$ in terms of the well-known generating function of the Catalan numbers, 
\begin{equation*} C(t) = \sum_{m = 0}^\infty C_m t^m = \frac{2}{1+\sqrt{1-4t}}.
\end{equation*}
\begin{lemma}\label{lemmagf}
If $p \geq 5$ is prime, then we have that
\[Z_p^{C}(t) = \sum_{m=0}^{\infty} M_p^{C}(m)t^{m} = (p-1)\left(pC(pt^2)-C(pt^2)^3 \frac{t^2}{1-t^2C(pt^2)^2}\right).\]
\end{lemma}
\begin{proof}
We use the identity \[\sum_{n = k}^\infty \frac{(2k+1)(2n)!}{(n-k)!(n+k+1)!} x^{n - k} = C(x)^{2k+1}.\] For $k = 0$, this is the usual generating function of the Catalan numbers. One can prove the general identity inductively by differentiating, multiplying by $x^{2k+2}$, and integrating. 

Given this identity, the formula for the generating function $Z_p^{C}(t)$  is a consequence of Birch's formula:
\begin{align*}
\sum_{m=0}^\infty \frac{M_p(m)}{p-1} t^m &= \sum_{m = 0}^\infty \frac{(2m)!}{m! (m+1)!} p^{m+1} t^{2m} - \sum_{m = 0}^\infty \sum_{k = 1}^m  \frac{(2k+1)(2m)!}{(m-k)!(m+k+1)!}p^{m-k}t^{2m} \\
&= pC(pt^2) - \sum_{k = 1}^\infty C(pt^2)^{2k + 1}t^{2k} \\
&= pC(pt^2) - C(pt^2)^3\cdot \frac{t^2}{1 - t^2C (pt^2)^2}.
\end{align*}
\end{proof}
Now we express $Z_p^H(t)$ in a convenient form in order to take advantage of known congruences for traces of Hecke operators. We rewrite the expression for $Z_p^H(t)$ as follows:
\begin{align*}
-\frac{Z_p^H(t)}{p-1} &= \sum_{m=0}^\infty \sum_{k=1}^m \frac{(2k+1)!(2m)!}{(m-k)!(m+k+1)!}p^{m-k}\Tr_{2k+2}(T_p) t^{2m} \\
&= \sum_{k=1}^\infty \Tr_{2k + 2}(T_p)C(pt^2)^{2k+1} t^{2k}.
\end{align*}
Thus, we have derived the following lemma. 
\begin{lemma}\label{heckegf}
If $p\geq 5$ is prime, then we have that
\begin{align*}
-\frac{Z_p^H(t)}{p-1} = \sum_{k=1}^\infty \Tr_{2k}(T_p)  C(pt^2)^{2k+1}t^{2k}.
\end{align*}
\end{lemma}

\begin{remark}
For $k \leq 5$, there are no cusp forms of weight $2k$, and so the contribution of $Z_p^H(t)$  to $Z_p(t)$ is $O(t^{12})$. Therefore, our closed form for $Z_p^C(t)$ gives precise values of $M_p(2m)$ for all $m \leq 5$.
\end{remark}

We now turn to congruence relations, dealing first with $M_p^C(t)$.
\begin{lemma}\label{congpoly} If $p \geq 5$ is prime and $m \geq 1$, then: 

\noindent(1) The following congruence modulo $p + 1$ holds: $$\frac{M_{p}^{C}(m)}{p-1} \equiv 0 \pmod{p+1}.$$ (2)  The following congruence modulo $p$ holds: \[M_p^C(2m) \equiv 1 \pmod{p}.\]
\end{lemma}

\begin{proof}
(1) To evaluate \[\sum_{m=0}^\infty M_p^C(2m)t^{2m} \pmod{p+1},\] we formally plug in $-1$ for $p$ in the closed form of the generating function. We obtain \[\frac{Z_{p}^{C}(t)}{p-1} \equiv  - C(-t^2) -C(-t^2)^3\frac{t^2}{1 - t^2 C(-t^2)^2} \equiv -1 \pmod{p+1}.\] Since all nontrivial terms of the generating function are 0 modulo $p+1$, the proof is complete.

\medskip
\noindent(2) Similarly, we  now formally plug in $0$ for $p$ in the closed form of the generating function. We see that $C(pt^2) \equiv 1 \pmod{p}$, so we have that 
\[
Z_p^C(t) \equiv \frac{t^2}{1-t^2} \equiv \sum_{m=1}^{\infty} t^{2m} \pmod{p}.
\]
So $M_p^{C}(2m) \equiv 1 \pmod{p}$. 
\end{proof}
We are now ready to prove Theorem \ref{congruencethm}. 
\begin{proof}For (1), if $\ell \in \left\{2, 3, 5\right\}$, $m \geq 1$, and $p \geq 5$ is a prime such that $p \equiv -1 \pmod{\ell}$, then we have by Lemma \ref{congpoly} that $${\frac{M_p^C(2m)}{p-1} \equiv 0 \pmod{\ell}}.$$  By a paper of Hatada \cite{Hatada}, for such $p$, all traces of Hecke operators satisfy $${\Tr_{2k} (T_p) \equiv 0 \pmod{\ell}}.$$ Since $C(pt^2)$ has integer coefficients, we see from Lemma \ref{heckegf} that $${\frac{M_p^H(2m)}{p-1} \equiv 0 \pmod{\ell}}.$$ Then the decomposition $M_p(2m) = M_p^H(2m) + M_p^C(2m)$ completes the proof of (1).

For (2), a result by Choie, Kohnen, and Ono \cite{CKO} shows that, under the hypotheses of the theorem, $$\Tr_m(T_p) \equiv 0 \pmod{p}.$$ Therefore we conclude, again by Lemmas  \ref{heckegf} and \ref{congpoly} and using the fact that $C(pt^2) \equiv 1 \pmod{p}$, that $$M_p(m) \equiv M_p^C (m) \equiv 1 \pmod{p}.$$
\end{proof}
 
\section{Examples} \label{sec4}
\begin{example}  We first illustrate Theorem \ref{TheoremZetaGF} for $A = \left\{0\right\}$ and $p = 5$. One can directly compute the power moments: $M_{5}(0)=20, M_{5}(2)=96, M_{5}(4)=936, M_{5}(6)=11496,$ and $M_{5}(8)=158856$. Then, Theorem \ref{TheoremZetaGF} implies that
$$\widehat{Z}_{5}(t)= 20-504t^{2} + 3516t^{4} - 6776t^{6} + 2304t^{8},$$ which in turn yields
\begin{align*}
Z_{5}(t) &= \frac{20-504t^{2} + 3516t^{4} - 6776t^{6} + 2304t^{8}}{1-30t^{2}+273t^{4}-820t^{6}+576t^{8}}\\
&= 20+96t^{2}+936t^{4}+ \cdots + 2212976684616t^{20}+ \cdots.
\end{align*}

One can easily confirm that $M_{5}(20)=2212976684616$ by direct computation. 
\end{example}
\begin{example}
We now show an example of Theorem \ref{TheoremTraceGF} when $A = \left\{0\right\}$ and $p = 5$:
\begin{align*}
\sum_{k =1}^\infty \Tr_{2k}(T_{5}) t^{2k} = 4830t^{12} + 52110t^{16} - 1025850t^{18} - 2377410t^{20} + 21640950t^{22} 
+ \cdots.
\end{align*}

\noindent The observation that the smallest power of $t$ appearing above is $12$ and that the coefficient of $t^{14}$ vanishes reflects that $S_{k}$ is trivial for $k \leq 10$ and $k=14$. It is well-known that $S_{k}$ is generated by a single form  for $k \in \left\{12,16,18,20,22\right\}$. We have $S_{12} = \langle \Delta\rangle$, $S_{16} = \langle E_4\Delta\rangle,$ $S_{18} = \langle E_6\Delta\rangle$, $S_{20} = \langle E_4^2 \Delta \rangle$, and $S_{22} = \langle E_4E_6\Delta\rangle$. Since each of these generators is a normalized eigenform, the trace of $T_5$ on $S_k$ for each of these spaces is its eigenvalue, which in turn is the coefficient of $q^5$ in the Fourier expansion. Therefore our generating function allows us to read off the coefficients of $q^5$ for $\Delta, E_4\Delta, \ldots, E_4E_6\Delta$. For instance, $\tau(5) = 4830$. 
\end{example}
\begin{example}
We demonstrate Theorem \ref{TheoremTraceGF} when $p \geq 5$ is a prime and $A = \mathbb{Z}/\ell\mathbb{Z}$, where $\ell$ is an odd prime such that $p \not\equiv 1 \pmod{\ell}$ and $p \neq \ell$. The relevant congruence subgroup is $\Gamma(\ell, 1) = \Gamma_1(\ell)$ and our trace operator formula is \begin{align*}\sum_{k = 2}^\infty \Tr\big(T_p | S_k \big(\Gamma_1(\ell)\big)\big) = (p+1) t^2 &- \frac{(\ell - 1) (1 + \delta_\ell (p, -1)) }{2} \cdot \frac{t^2}{1-t^2} \\ &-  \frac{(\ell-1)(1 - \delta_\ell (p, -1)) }{2} \cdot \frac{t^3}{1 - t^2} \\ &- \frac{(\ell - 1)t^2}{(p - 1)(1 + pt^2)}\cdot Z_p\left(\mathbb{Z}/\ell\mathbb{Z}; \frac{t}{1 + pt^2}\right).\end{align*}
\end{example}
\begin{example}
As another example of Theorem \ref{TheoremTraceGF}, take a prime  $p \geq 5$  and let $A = \mathbb{Z}/\ell\mathbb{Z} \times \mathbb{Z}/\ell\mathbb{Z}$ for an odd prime $\ell$ such that $p \equiv 1 \pmod{\ell}$. The relevant congruence subgroup is $\Gamma(\ell, \ell)$, which is isomorphic to $\Gamma(\ell)$. The theorem yields \begin{align*}\sum_{k = 2}^\infty \Tr\big(T_p | S_k \big(\Gamma(\ell)\big)\big) t^k = (p + 1) \cdot t^2 &- \frac{(t^2 + t^3)(\ell^2 - 1)}{2(1 - t^2)} \\ &- \frac{(\ell^3 - \ell)t^2}{(p-1)(1 + pt^2)} \cdot Z_p\left(\mathbb{Z}/\ell\mathbb{Z} \times \mathbb{Z} /\ell\mathbb{Z}; \frac{t}{1 + pt^2}\right). \end{align*}
\end{example}


\begin{document}
\maketitle

\begin{abstract}
Seminal works by Birch and Ihara gave formulas for the $m$th power moments of the traces of Frobenius endomorphisms of elliptic curves over $\mathbb{F}_{p}$ for primes $p \geq 5$. Recent works by Kaplan and Petrow generalized these results to the setting of elliptic curves that contain a subgroup isomorphic to a fixed finite abelian group $A$. We revisit these formulas and determine a simple expression for the zeta function $Z_p(A; t)$, the generating function for these $m$th power moments. In particular, we find that
\[
Z_p(A;t) = \frac{\widehat{Z}_p(A; t)}{\hspace*{-0.3cm}\displaystyle \prod_{a \in \Frob_p(A)}\hspace*{-0.3cm}(1 - at)},\]
where $\Frob_p(A) \coloneqq \{ a \, \colon -2\sqrt{p} \leq a \leq 2\sqrt{p}\, \text{ and }  a \equiv p+1 \pmod{|A|}\}$, and $\widehat{Z}_p(A;t)$ is an easily computed polynomial that is determined by the first $\Big\lceil\frac{2\lfloor 2\sqrt{p}\rfloor}{|A|}\Big\rceil$ power moments.
These rational zeta functions have two natural applications. We find rational generating functions in weight aspect for traces of Hecke operators on $S_k(\Gamma)$ for various congruence subgroups $\Gamma$. We also prove congruence relations for power moments by making use of known congruences for traces of Hecke operators.

\end{abstract}






\section*{Acknowledgements}

\thanks{The authors wish to thank Professor Ken Ono and Professor Larry Rolen for their invaluable guidance and suggestions. They would also like to thank Professor Nathan Kaplan and Professor Ian Petrow for taking the time to read our paper and provide helpful feedback. They also thank Emory University, the Asa Griggs Candler Fund, and NSF grant DMS-1557960.}
\bibliographystyle{abbrv}
\bibliography{biblio}

\begin{thebibliography}{10}

\bibitem{Birch}
B.~J. Birch.
\newblock How the number of points of an elliptic curve over a fixed prime
  field varies.
\newblock {\em J. London Math. Soc.}, 43:57--60, 1968.

\bibitem{CKO}
Y.~Choie, W.~Kohnen, and K.~Ono.
\newblock Linear relations between modular form coefficients and non-ordinary
  primes.
\newblock {\em Bull. London Math. Soc.}, 37(3):335--341, 2005.

\bibitem{STI}
L.~Clozel, M.~Harris, and R.~Taylor.
\newblock Automorphy for some $\ell$-adic lifts of automorphic $\mod \ell$
  galois representations.
\newblock {\em Publications Math\'ematiques de l'IH\'ES}, 108(1):1--181, 2008.

\bibitem{STII}
L.~Clozel, M.~Harris, and R.~Taylor.
\newblock Automorphy for some {$l$}-adic lifts of automorphic mod {$l$}
  {G}alois representations.
\newblock {\em Publ. Math. Inst. Hautes \'Etudes Sci.}, (108):1--181, 2008.
\newblock With Appendix A, summarizing unpublished work of Russ Mann, and
  Appendix B by Marie-France Vign\'eras.

\bibitem{Cox}
D.~A. Cox.
\newblock {\em Primes of the Form $x^2 + ny^2$}.
\newblock John Wiley Sons, \& Inc, 1989.

\bibitem{Deuring}
M.~Deuring.
\newblock Die typen der multiplikatorenringe elliptischer functionenk{\"o}rper.
\newblock {\em Abh{.} Math{.} Sem{.} Hansischen Univ{.}}, 14:197--272, 1941.

\bibitem{DS}
F.~Diamond and J.~Shurman.
\newblock {\em A First Course in Modular Forms}, volume 228.
\newblock Springer-Verlag New York, New York, 2005.

\bibitem{Hatada}
K.~Hatada.
\newblock Congruences for eigenvalues of {H}ecke operators on {${\rm
  SL}_{2}({\bf Z})$}.
\newblock {\em Manuscripta Math.}, 34(2-3):305--326, 1981.

\bibitem{Ihara}
Y.~Ihara.
\newblock Hecke {P}olynomials as congruence {$\zeta $} functions in elliptic
  modular case.
\newblock {\em Ann. of Math. (2)}, 85:267--295, 1967.

\bibitem{KaplanPetrow}
N.~Kaplan and I.~Petrow.
\newblock Elliptic curves over a finite field and the trace formula.
\newblock {\em Proc. Lond. Math. Soc. (3)}, 115(6):1317--1372, 2017.

\bibitem{Kowalski}
E.~Kowalski.
\newblock Analytic problems for elliptic curves.
\newblock {\em J. Ramanujan Math. Soc.}, 21(1):19--114, 2006.

\bibitem{Ono}
K.~Ono.
\newblock {\em The Web of Modularity: Arithmetic of the Coefficients of Modular
  Forms and q-series}, volume 102.
\newblock CBMS Regional Conference Series in Mathematics{,} Amer. Math. Soc.,
  Providence{,} RI, 2004.

\bibitem{STIII}
R.~Taylor.
\newblock Automorphy for some $\ell$-adic lifts of automorphic mod $\ell$
  galois representations. {II}.
\newblock {\em Publications Math\'ematiques de l'IH\'ES}, 108(1):183--239,
  2008.

\end{thebibliography}

\end{document}